\documentclass[A4,11pt]{article}
\usepackage{amsmath}
\usepackage{amssymb}
\usepackage{amscd}
\newtheorem{theorem}{Theorem}[section]
\newtheorem{proposition}[theorem]{Proposition}
\newtheorem{corollary}[theorem]{Corollary}
\newtheorem{lemma}[theorem]{Lemma}

\newtheorem{remarks}[theorem]{Remarks}
\newenvironment{proof}{\begin{trivlist} 
\item[] {\bf Proof.}}
{\hfill $\square$\end{trivlist}}

\newcommand{\DSH}{{\rm DSH}}

\newcommand{\Ccal}{{\cal C}}

\renewcommand{\P}{\mathbb{P}}
\newcommand{\C}{\mathbb{C}}

\newcommand{\R}{\mathbb{R}}

\newcommand{\Z}{\mathbb{Z}}
\newcommand{\id}{{\rm id}}
\newcommand{\h}{{\rm h}}

\renewcommand{\H}{{\cal H}}

\newcommand{\K}{{\cal K}}

\newcommand{\dc}{{\rm d^c}}
\newcommand{\ddc}{{\rm dd^c}}
\renewcommand{\d}{{\rm d}}

\newcommand{\Lone}{{{\rm L}^1}}
\newcommand{\Ltwo}{{{\rm L}^2}}

\newcommand{\D}{{\cal D}}

\renewcommand{\O}{{\rm O}}

\newcommand{\End}{{\rm End}}

\title{Green currents for holomorphic automorphisms of
compact K\"ahler manifolds}
\author{Tien-Cuong Dinh and Nessim Sibony}
\begin{document}
\maketitle
\begin{abstract} Let $f$ be a holomorphic automorphism of a compact 
K\"ahler manifold $(X,\omega)$ of dimension $k\geq 2$. 
We study the convex cones of positive 
closed $(p,p)$-currents $T_p$, which 
satisfy a functional relation
$$f^* T_p=\lambda T_p,\ \ \lambda>1,$$
and some regularity condition (PB).
Under appropriate assumptions on dynamical degrees we introduce closed finite 
dimensional cones, not reduced to zero, of such currents. 
In particular, when the topological entropy $\h(f)$ of $f$ is positive, 
then for some $m\geq 1$, there is a positive closed $(m,m)$-current $T_m$
which satisfies the relation
$$f^* T_m=\exp(\h(f)) T_m.$$
Moreover,
every quasi-p.s.h. function is integrable with respect to 
the trace measure of $T_m$. 
When the 
dynamical degrees of $f$ are all distinct, we construct 
an invariant measure $\mu$ as an intersection 
of closed currents. We show that this measure is mixing and
gives no mass to pluripolar sets and 
to sets of small Hausdorff dimension. 
\end{abstract}
\section{Introduction}
Let $(X,\omega)$ be a compact K\"ahler manifold of dimension $k$. 
Let $f$ be a holomorphic automorphism
of $X$. 
Our purpose is to introduce 
invariant positive closed currents for $f$ and to use them in order 
to construct dynamically interesting invariant measures.
\par
More precisely we want to construct positive closed $(s,s)$-currents $T$ 
satisfying a functional equation 
$$f^*T=\lambda T,\ \ \lambda>1,$$
and some additional regularity properties.
Very likely, the currents 
$T$ will describe the distribution of invariant manifolds 
of codimension $s$ corresponding to the smallest Lyapounov exponents. 
\par
Let $d_p$ denote the dynamical degree of order $p$ of $f$.
It describes the growth under iteration of the volume of $p$-dimensional
manifolds.
When $d_1>1$, 
it is natural to introduce first a positive closed $(1,1)$-current as
$$T_1:=\lim_{n\rightarrow\infty} \frac{(f^n)^*\omega}{d_1^n}.$$
Unfortunately, this limit does not always exist, as holomorphic automorphisms
of tori show.
\par
However, for H\'enon maps in $\C^2$ or for algebraically stable meromorphic 
self-maps 
of $\P^k$ the limit exists and was studied extensively. See the introduction
of \cite{BedfordLyubichSmillie4, Sibony2} for historical comments. The limit
$T_1$ is called {\it the Green current}. For holomorphic endomorphisms 
of $\P^k$, the self-intersection $T_1^k$ is well defined and gives an 
invariant measure.
\par
In the context of polynomial automorphisms $f$ of $\C^2$ having positive entropy,
the second author has defined, using the intersection of 
the Green currents of $f$ and $f^{-1}$,
an invariant probability measure which turned out to be dynamically 
interesting. It, as well as
the Green current, was extensively studied by 
Bedford-Lyubich-Smillie \cite{BedfordSmillie, 
BedfordLyubichSmillie4} and by Forn\ae ss and the second author 
\cite{FornaessSibony1, FornaessSibony2}.
\par
Cantat \cite{Cantat} 
has adapted some of the constructions from the study of H\'enon maps
to the context of automorphisms of K3 surfaces (see also \cite{Guedj1}). 
The striking fact in this 
case is the existence of automorphisms of positive entropy for some K3 
surfaces \cite{Mazur,Cantat,McMullen}.
\par
Given a holomorphic automorphism or, more generally, a bimeromorphic
self-map of a compact K\"ahler manifold of 
dimension $k\geq 3$,
constructing interesting invariant measures is a delicate process.
The case of  
some classes of polynomial automorphisms of $\C^k$ is considered in 
\cite{Sibony2,GuedjSibony}. 
The same difficulty arises for meromorphic self-maps whose
topological degree is not the largest of the dynamical degrees,
see \cite{FavreGuedj,Guedj1, DinhSibony1}.
In this case, the cohomology class $[T_1]^k$ is zero in the cohomology group 
$\H^{k,k}(X,\R)$. Hence we cannot construct an invariant measure as 
the self-intersection of the Green current $T_1$.
\par
In this paper, we focus on the question of constructing invariant currents and 
invariant measures for automorphisms of compact K\"ahler manifolds. We 
consider first automorphisms to avoid some technicalities. Some of the results
can be extended to holomorphic maps, to non singular correspondences and to 
birational maps \cite{DinhSibony7}.
\par 
When $d_1>1$, we consider the convex
cone $\Gamma_1$ of positive closed $(1,1)$-currents $T_1$ with bounded 
potentials such that $f^*T_1=d_1T_1$. We will show that this cone is closed, finite 
dimensional and not reduced to zero. Moreover, every current in $\Gamma_1$
has a H\"older continuous potential. An element of $\Gamma_1$ is called
{\it a Green current of bidegree $(1,1)$} (see \cite{Cantat, Guedj1} for the
case of surfaces).
\par
To construct invariant currents of higher bidegree, we start with a positive
closed current $T$ of bidegree $(s,s)$ satisfying 
$$f^*T=\lambda_T T, \ \lambda_T>0.$$ 
We introduce the cone $\Gamma(T)$ of 
positive closed $(s+1,s+1)$-currents $T\wedge S$ 
such that $T\wedge f^*S=\lambda_1(T)T\wedge S$
where $S$ is a closed (not necessarily positive)
$(1,1)$-current with a continuous potential and 
$$\lambda_1(T):=\lim_{n\rightarrow\infty} 
\left(\int T\wedge(f^n)^*\omega\wedge 
\omega^{k-s-1}\right)^{1/n}.$$
The number $\lambda_1(T)$ appears as a dynamical degree with respect to $T$.
When $\lambda_1(T)>1$, we show that $\Gamma(T)$ 
is a closed, finite dimensional, convex
cone, not reduced to zero. This requires the introduction of
cohomology groups with respect to $T$. We also prove 
that $\Gamma(T)$ is unchanged if we consider 
only currents $S$ with H\"older
continuous potentials. In this approach, 
we use an inductive procedure, and at the step $s$ 
we replace the complex manifold $X$ by the $(s,s)$-current $T$. 
\par
When $d_1$ is strictly larger than the other dynamical degrees of $f$, using 
$f$ and $f^{-1}$, we can
construct invariant positive currents 
as an intersection of 
closed $(1,1)$-currents. 
We first construct the current $T_1$ such that $f^*T_1=d_1T_1$ as above. We have 
$f_*T_1=d_1^{-1} T_1$. We then construct $(1,1)$-currents $S_i$, $1\leq i\leq k-1$,  
with H\"older continuous potentials such that the currents $T_1\wedge S_1 \wedge
\ldots\wedge S_i$ are
positive, closed and satisfy
$$f_*(T_1\wedge S_1\wedge\ldots\wedge S_i) = c_i T_1\wedge S_1
\wedge\ldots\wedge S_i \mbox{ for some }
c_i>0.$$
In this case, we can show, at each step,  that for the automorphism $f^{-1}$, 
the dynamical degree  
$\lambda_1(T_1\wedge S_1\wedge \ldots\wedge S_i)$, with respect to the current 
$T_1\wedge S_1\wedge \ldots\wedge S_i$, is 
strictly larger than 1.   
The measure $\mu:=T_k$ is invariant, mixing and gives
no mass to pluripolar sets and to sets of small Hausdorff dimension. 
\par
In the general case, the method breaks down because the hypothesis on 
$\lambda_1(T)$
is not easy to check. Therefore, we introduce a second
method for constructing directly an invariant current of bidegree
$(s,s)$, under the assumption $d_s>d_{s-1}$. The current we
construct is PB and satisfies 
\begin{equation}
f^*T  =  d_s T.
\end{equation} 
Note that if a current $T$ is PB, it is weakly PB, i.e.
q.p.s.h. functions are integrable with 
respect to the trace measure $T\wedge \omega^{k-s}$ (see Section 2.1 for details).
The method uses a {\it $\ddc$-resolution} as
already used by the authors in various contexts
\cite{DinhSibony2, DinhSibony3, Dinh}. 
We have to use complex analysis, i.e. estimates for the 
solutions of the $\ddc$-equation (Proposition 2.1). 
This permits one to control the growth of $(f^n)_*\Phi$, 
where $\Phi$ is a test form.
The delicate point in the construction is to obtain a PB current.
See also \cite{DinhSibony9} for another new method.

The method permits one to prove that $T$ is almost extremal, i.e. $T$ belongs to a
finite dimensional extremal face of the cone of positive closed currents.
\par
To construct an invariant measure $\mu$, we assume that the
dynamical degrees of the automorphism $f$ are distinct. Then the 
Khovanskii-Tessier-Gromov concavity theorem implies the existence of
an $m$ such that
$$1<d_1<\cdots<d_m>d_{m+1}>\cdots>d_k=1.$$
The measure $\mu$ is then obtained using the first method for $f^{-1}$ but
starting with $T_m$. Therefore, $\mu$ is a wedge product of $T_m$ with
closed $(1,1)$-currents with H\"older continuous potentials. We can choose
$T_m$ so that
$\mu$ is mixing and has positive Hausdorff dimension.
Observe that according to Yomdin-Gromov \cite{Yomdin, Gromov1, DinhSibony4, 
DinhSibony6} the topological
entropy $\h(f)$ of $f$ is equal to $\max \log d_i$. 
Hence, the PB current $T_m$ satisfies
$f^* T_m=\exp(\h(f)) T_m$. 

The classes of PC and PB currents, which we introduce, are of interest since 
they allow to define the product of currents of higher bidegree. We will 
come back to this subject in a furture work.

Many questions concerning Green 
currents and the measure $\mu$ are not studied here for instance: 
distribution of periodic points with respect to $\mu$, entropy of $\mu$, 
approximation of the Green currents by stable leaves... 
If $X$ is projective, the first author has proved 
that $T_m$ is weakly laminar \cite{Dinh2}.

The classification of manifolds of dimension $\geq 3$ with automorphisms
of positive entropy is also an interesting problem. In dimension 2, many
examples are given in \cite{Mazur, Cantat, McMullen}. Mazur's examples
can be extended to dimension $\geq 3$. 
\\

\noindent
{\bf Mazur's examples.}  
Let $X$ be a smooth hypersurface of degree 2 of $\P^{1,k+1}:=
\P^1\times \cdots
\times \P^1$ ($k+1$ factors). Let $\pi_i$, $1\leq
i\leq k+1$, denote the $k+1$ projections of $X$ on the product
of $k$ factors of $\P^{1,k+1}$. Assume that all $\pi_i$ 
are finite. Then, each fiber of
$\pi_i$ contains exactly 2 points $z$, $z'$. We can define
an involution $\tau_i$ by $z\mapsto z'$. The group generated by
$\tau_i$ contains dynamically interesting 
automorphisms.
\\

One can construct other examples by taking
products of manifolds or the projectivization of their tangent bundle. 
The dynamics of these examples is however non trivial. It is used in
\cite{Dinh2} to get information on laminarity of currents. 
Examples on tori explain some of the difficulties
that we have to overcome in the general case.
Our results extend to non singular correspondences (see \cite[p.337]{Margulis} 
and \cite{ClozelUllmo, Voisin, Dinh, DinhSibony3} for definitions and examples).
\section{Currents and groups of cohomology}
A holomorphic automorphism $f$ of $X$ 
induces an invertible linear self-map on
groups of cohomology. We will use this action of $f$ 
in order to construct
invariant currents.  
We introduce classes of currents with some regularity properties.
\par
We will write
$u_n\simeq v_n$ for
$\lim u_n/v_n=1$.
The notation $u_n\lesssim v_n$ means 
$\limsup |u_n/v_n|<+\infty$, the notation $u_n\sim v_n$ 
means $u_n\lesssim v_n$ and $v_n\lesssim u_n$, with the convention that 
$0\simeq 0$, $0\lesssim 0$ and $0\sim 0$. For $(x,y)$ and $(x',y')$
in $\R^2$, we write $(x,y)\leq (x',y')$ if either $x<x'$, or if $x=x'$
and $y\leq y'$. The sign $\|\ \|$ denotes either
the the mass of currents, the norm of a vector, or of a linear operator. 
The sign $[\ ]$ 
denotes a cohomology class of a closed 
current.
\\
\ 
\\
{\bf 2.1 PC, PB and weakly PB currents}
\\
\ 
\\
We refer to the 
survey by Demailly \cite{Demailly1} for the basics 
on currents in complex analysis.
Demailly's survey on Hodge theory \cite{Demailly2} is also useful.
We however recall a few facts. 
\par 
When $(X,\omega)$ 
is a compact K\"ahler manifold of dimension 
$k$, a {\it current $T$ of bidegree $(s,s)$} is a continuous linear form 
on $\D^{k-s,k-s}(X)$ -- the space of smooth forms of bidegree $(k-s,k-s)$. In a 
coordinate chart, $T$ can be expressed as a differential 
$(s,s)$-form with distribution 
coefficients.
A $(k-s,k-s)$-form $\Phi$ is {\it weakly positive} if at every point $z\in X$
$$\Phi\wedge i\alpha_1\wedge \overline\alpha_1\wedge \ldots \wedge i\alpha_s\wedge 
\overline \alpha_s\geq 0$$
for every $(1,0)$-form 
$\alpha_j$
of $X$. The current $T$ is called {\it (strongly) positive} if 
$\langle T, \Phi\rangle\geq 0$  for every weakly positive test form $\Phi$. 
The space of currents is given the weak topology of currents. 
In particular, positive 
currents are currents of order zero. The {\it trace measure} $\sigma_T$ 
associated to a positive current $T$ is defined as 
$\sigma_T:=\frac{1}{(k-s)!}T\wedge \omega^{k-s}$. 
The measure $\sigma_T$ is positive 
and the coefficients of $T$ in a chart are measures 
which are dominated by $c\sigma_T$, $c>0$. We will denote by $\|T\|$ {\it the mass} 
$\int_X T\wedge \omega^{k-s}$ of $T$.
\par 
The calculus on differential forms extends to currents, except for the 
pullback by a holomorphic map, which is not a submersion. It is always delicate
to define the wedge product of two currents. However when $u\in\Lone(\sigma_T)$,
for example if $u$ is bounded, one can define $\ddc u\wedge T:=\ddc (uT)$.
The continuity properties of this operator depend on the properties of $u$ 
\cite{Demailly1}. Recall that 
$\d=\partial+\overline\partial$, $\dc=\frac{1}{2i\pi}
(\partial -\overline\partial)$ and  
that $\ddc=\frac{i}{\pi}\partial\overline\partial$ 
is a real operator.
\par
Let $T$ be a positive closed $(s,s)$-current, and $S$ be a closed 
$(1,1)$-current not necessarily positive.
Since $X$ is K\"ahler, by Hodge theory,  
we can write $S=\alpha +\ddc u$, where $\alpha$ is a smooth form cohomologous
to $S$ and $u$ is a $(0,0)$-current. We say that $u$ is {\it a potential} 
of $S$. Observe that two potentials of $S$ differ by a smooth function.
When $u$ is a
$\nu$-H\"older continuous (resp. continuous or bounded) function, 
we say that $S$ has {\it a $\nu$-H\"older continuous}
(resp. {\it a continuous or bounded}) {\it potential}.
It is clear that this is independent of the choice of $\alpha$.
For a current $S$ with a continuous potential, we can define $T\wedge S$ by 
$T\wedge S:=T\wedge\alpha+\ddc(uT)$. When $S$ is positive, we can choose $u$ 
upper semi-continuous. In this case, if $u$ is bounded, 
we can define $T\wedge S$ 
in the same way. 
\par

A real $(p,p)$-current $\Phi$ on $X$ is called {\it DSH} if $\Phi=\Phi_1-\Phi_2$ 
and  
$\ddc\Phi_i=\Omega^+_i-\Omega^-_i$ where 
$\Phi_i$ are negative currents, and $\Omega^\pm_i$ are positive closed currents. 
We define the DSH-norm by 
$$\|\Phi\|_\DSH:=\inf\big\{\|\Phi_1\|+\|\Phi_2\| + \|\Omega_1^+\|
+\|\Omega_2^+\|,\  
\Phi_i,\ \Omega_i^\pm \mbox{ as above}\big\}.$$
Observe that $\|\Omega^+_i\|=\|\Omega^-_i\|$, and we can choose $\Phi_i$ and 
$\Omega_i^\pm$ such
that $\|\Phi\|_\DSH=\|\Phi_1\|+\|\Phi_2\| + \|\Omega_1^+\|
+\|\Omega_2\|$. 
Denote by $\DSH^p(X)$ the space of 
DSH $(p,p)$-currents. This is our space of test currents.
A current $\Phi$ is in $\DSH^0(X)$ if and only if it is a Difference of q.p.S.H. 
functions.

Recall that an $\Lone$ function  $\varphi:X\rightarrow \R\cup\{-\infty\}$ 
is {\it quasi-plurisubharmonic} ({\it q.p.s.h.}
for short) if 
it is upper semi-continuous and $\ddc\varphi\geq -c\omega$, $c>0$, in the sense of
currents. A set $E\subset X$ is {\it pluripolar} if
it is contained in the pole set 
$\{\varphi=-\infty\}$ of a q.p.s.h. function $\varphi$. 

A topology on $\DSH^p(X)$ is defined in the following way: $\Phi^{(n)}\rightarrow 
\Phi$ in  $\DSH^p(X)$ if $\Phi^{(n)}\rightarrow 
\Phi$ weakly and if $(\|\Phi^{(n)}\|_\DSH)$ is bounded.

Let $T$ be a current of bidegree $(s,s)$ and of zero order. We say that 
$T$ is {\it PC} if it can be extended to a continuous linear form on 
$\DSH^{k-s}(X)$, and we write $\langle T,\Phi\rangle$ for the value of this linear 
form on $\Phi\in \DSH^{k-s}(X)$. In \cite{DinhSibony6}, we proved that every current
in $\DSH^{k-s}(X)$ can be approximated by smooth forms. Then, the extension of $T$
is unique. Moreover, when $\Phi$ is continuous  $\langle T,\Phi\rangle$ coincides 
with the usual integral \cite[Prop. 4.6]{DinhSibony6}.

We say that $T$ is {\it PB} if there exists a constant $c_T>0$ such that
$$|\langle T,\Phi\rangle| \leq c_T\|\Phi\|_\DSH
\ \mbox{ for every DSH continuous } (k-s,k-s)\mbox{-form }\Phi.$$
The current $T$ is called {\it weakly PB} 
if $$|\langle T,\varphi\omega^{k-s}\rangle|<+\infty
\ \mbox{ for every q.p.s.h. function } \varphi \mbox{ on  }X.$$
 
Observe that $T$ is weakly PB if and only if the measure $T\wedge \omega^{k-s}$ is
weakly PB. If a positive 
current $T$ is (weakly) PB,
then every positive current $T'$ such that $T'\leq T$, is also (weakly) PB. 
We can show that if $T$ is a positive weakly PB current, then
$|\langle T,\varphi\omega^{k-s}\rangle|\leq c_T(1+\|\varphi\|_\Lone)$ for some
constant $c_T>0$ and for $\varphi$ q.p.s.h. such that 
$\ddc\varphi\geq -\omega$ (see \cite{DinhSibony3} or Proposition 2.2).

In \cite{DinhSibony2}, we showed that
a positive measure $\mu$ on a Riemann surface admits locally 
a Bounded Potential if and only 
if, locally p.s.h. functions are $\mu$-integrable.
This justifies the term PB (see also Proposition 2.2).
\par
The following result is useful in constructing PC and PB currents.
It supplies the fact that when $\Omega^-$ is a smooth form cohomologous to a positive
closed current $\Omega^+$, one cannot find, in general, a negative form $\Phi$
solving $\ddc\Phi=\Omega^+-\Omega^-$. A counter-example can be founding in 
\cite{BostGilletSoule}.  
\begin{proposition} There exists a constant $A>0$ so that for every pair of 
positive
closed $(k-s+1,k-s+1)$-currents $\Omega^\pm$ on $X$ with $[\Omega^+]=[\Omega^-]$, 
there are $\Lone$ negative $(k-s,k-s)$-forms $\Phi^\pm$ such that 
$$\ddc\Phi^+-\ddc\Phi^-=\Omega^+-\Omega^-
\ \mbox{ and }\ \|\Phi^\pm\|_\DSH\leq A\|\Omega^+\|.$$ 
Moreover, the DSH currents $\Phi^\pm$ depend continuously on $\Omega^\pm$. 
If $\Omega^\pm$ are continuous, then $\Phi^\pm$ are continuous.
\end{proposition}
\begin{proof} 
Analogous problems are considered in \cite{DinhSibony6} 
and other aspects of the following computation are detailed there.
By Hodge theory \cite{GriffithsHarris}, we have
$$\H^{k,k}(X\times X,\C)\simeq \sum_{p+p'=k\atop q+q'=k} \H^{p,q}(X,\C) \otimes_\C
\H^{p',q'}(X,\C).$$
Hence, 
if $\Delta$ is the diagonal of $X\times X$,
there exists a smooth real $(k,k)$-form $\alpha(x,y)$ on $X\times X$,
cohomologous to $[\Delta]$ with
$\d_x\alpha=\d_y\alpha=0$.
Following Bost-Gillet-Soul\'e 
\cite{GilletSoule, BostGilletSoule}, 
one can construct a $(k-1,k-1)$-form $K(x,y)$ on 
$X\times X$ such that 
$\ddc K=[\Delta]-\alpha$. We recall the construction.
\par
Let  $\pi:\widehat{X\times X}\longrightarrow 
X\times X$ be the blow-up of $X\times X$ along $\Delta$. 
It follows from a theorem of Blanchard \cite{Blanchard} 
that $\widehat{X\times X}$ is a K\"ahler manifold. 
Let $\widehat\Delta:=\pi^{-1}(\Delta)$
be the exceptional hypersuface. 
Choose a negative q.p.s.h. function $\widehat\varphi$
on $\widehat{X\times X}$ such that
$\gamma:=-\ddc\widehat\varphi +[\widehat\Delta]$ is a smooth form 
cohomologous to $[\widehat\Delta]$. 
We can choose \cite[1.3.6]{GilletSoule} a smooth closed real $(k-1,k-1)$-form
$\eta$ on $\widehat {X\times X}$ such that $\pi^*\alpha$ is cohomologous to
$[\widehat\Delta]\wedge \eta$, hence to $\gamma\wedge\eta$.  
It follows that there is a smooth real $(k-1,k-1)$-form $\beta$ such that
$\ddc\beta=-\gamma\wedge \eta +\pi^*\alpha$. 
Define 
$$K(x,y):=\pi_*(\widehat\varphi \eta -\beta).$$
We have
$$\ddc K(x,y)=\pi_*([\widehat\Delta]\wedge\eta - \pi^*\alpha)
= \pi_*([\widehat\Delta]\wedge\eta) - \alpha.$$
The  $(k,k)$-current $\pi_*([\widehat\Delta]\wedge\eta)$ is 
closed, of order zero and supported on $\Delta$. Hence, it is 
a multiple of $[\Delta]$. Moreover, $\pi_*([\widehat\Delta]\wedge\eta)$
is cohomologous to $[\alpha]=[\Delta]$. 
It follows that $\pi_*([\widehat\Delta]\wedge\eta)
=[\Delta]$ and $\ddc K=[\Delta]-\alpha$. 
\par
Since the forms $\eta$ and $\beta$ are smooth, we can write
$\eta:=\eta^+-\eta^-$ and $\beta:=\beta^+-\beta^-$ with positive closed smooth
forms $\eta^\pm$ and negative smooth forms $\beta^\pm$. Define
$K^\pm:=\pi_*(\widehat\varphi \eta^\pm+\beta^\mp)$. These forms are negative and
we have $K=K^+-K^-$. Moreover, there exist constants $c^\pm$ with $c^+-c^-=1$
and closed real $(k,k)$-forms $\Theta^\pm$ on $X\times X$
such that $\ddc K^\pm=\Theta^\pm+c^\pm[\Delta]$.  
Let $|x-y|$ denote the distance between two points $x$ and $y$ of $X$ with
respect to the K\"ahler metric on $X$. 
One can check that $K^\pm$, $\Theta^\pm$ are smooth 
on $X\times X\setminus\Delta$ and
$$K^\pm(x,y)=\O\big(|x-y|^{2-2k}\log|x-y|\big), \ \ \  
\Theta^\pm=\O\big(|x-y|^{2-2k}\big)$$ 
when  $(x,y)\rightarrow \Delta$.
This allows one to define
$$\Phi^+(x):=\int_{y\in X}   K^+(x,y)\wedge \Omega^+(y)
+\int_{y\in X}   K^-(x,y)\wedge \Omega^-(y)$$
$$\Phi^-(x):=\int_{y\in X}   K^+(x,y)\wedge \Omega^-(y)
+\int_{y\in X}   K^-(x,y)\wedge \Omega^+(y).$$
Since $\ddc K=[\Delta]-\alpha$, 
$\d_y\alpha=0$, and $\Omega^+ -\Omega^-$ is exact, we have
\begin{eqnarray*}
\lefteqn{\ddc\Phi^+(x) -\ddc\Phi^-(x)  =  \int_{y\in X} (\ddc)_x 
K(x,y)\wedge \big(\Omega^+(y) -\Omega^-(y)\big)}\\
&\hspace{1cm} = & \int_{y\in X} \ddc 
K(x,y)\wedge \big(\Omega^+(y) -\Omega^-(y)\big)\\
&\hspace{1cm} = & \Omega^+(x) -\Omega^-(x)-\int_{y\in X} \alpha(x,y)\wedge 
\big(\Omega^+(y) -\Omega^-(y)\big)\\
&\hspace{1cm} = & \Omega^+(x) -\Omega^-(x).
\end{eqnarray*}
\par
The description of the singularities of $K^\pm$ 
implies that $\|\Phi^\pm\|_\Lone\lesssim \|\Omega^+\|$, that $\Phi^\pm$ depend 
continuously on $\Omega^\pm$
and that  $\Phi^\pm$ are continuous when 
$\Omega^\pm$ are continuous.
We can also write $\Theta^\pm$ as differences of positive closed forms, smooth
on $X\times X\setminus\Delta$, with singularities of order $\O(|x-y|^{2-2k})$. 
It follows that $\|\Phi^\pm\|_\DSH\lesssim \|\Omega^+\|$.
\end{proof}
\begin{proposition} 
Let $T$ be a positive $(s,s)$-current on $X$. If $T$ is PC, then $T$ is PB.
If $T$ is PB, then it is weakly PB. 
A positive measure $\mu$ on $X$ is PB if and only if it is weakly PB.
\end{proposition}
\begin{proof}
Let $\Phi_n$ be continuous forms with $\|\Phi_n\|_\DSH=1$.
If $c_n:=|\langle T,\Phi_n\rangle|$ tend to $+\infty$, then $c_n^{-1}\Phi_n$ converge
to $0$ in $\DSH^{k-s}(X)$, and $|\langle T,c_n^{-1}\Phi_n\rangle|$ converge to 1.
Hence, $T$ is not PC.

Assume now that $T$ is PB. 
Let $\varphi$ be a strictly negative q.p.s.h. function on $X$
such that $\ddc\varphi\geq -\omega$.
We have to show that $\langle T,\varphi\omega^{k-s}
\rangle >-\infty$. 
Following a theorem of Demailly \cite{Demailly3}, there exists a decreasing 
sequence of smooth negative functions $(\varphi_n)$ with limit $\varphi$, 
which satisfy $\ddc\varphi_n \geq -c_X \omega$ where $c_X>0$ is a constant.
We then have
$$|\langle T,\varphi\omega^{k-s}\rangle| = \lim |\langle T,\varphi_n \omega^{k-s}
\rangle| \lesssim \limsup\|\varphi_n\omega^{k-s}\|_\DSH 
\lesssim \|\varphi\|_\Lone +1.$$   
Hence $T$ is weakly PB.
\par
Now, assume that $\mu$ is a weakly 
PB probability measure. We show that it is PB. 
If not, there would exist continuous functions $\varphi_n$ such that
$\|\varphi_n\|_\DSH=1$ and  $\langle \mu, \varphi_n\rangle \geq n^3$.
We can write $\ddc\varphi_n = \Omega_n^+-\Omega_n^-$ where $\Omega_n^\pm$ are 
positive closed $(1,1)$-currents such that $\|\Omega_n^\pm\|\leq 1$. 
Since the set of such currents is compact, there exist
smooth forms $\Omega_n$, cohomologous to $\Omega_n^\pm$, 
such that $\Omega_n\leq c\omega$,
where $c>0$ is a constant independent of $n$.  
There exist q.p.s.h. functions $\varphi^\pm_n$ satisfying 
$\ddc\varphi^\pm_n =\Omega^\pm_n
-\Omega_n$ and $\max_X\varphi^\pm_n=0$. The family of q.p.s.h. functions 
$\psi$ such that $\max_X\psi=0$ and $\ddc\psi\geq -c\omega$ is compact in 
$\Lone(X)$. Hence there is a constant $A>0$ such that 
$\|\varphi^\pm_n\|_\DSH\leq A$. Define $c_n:=\varphi_n-\varphi_n^++\varphi_n^-$.
We have $\ddc c_n=0$. Hence $c_n$ is a constant and $|c_n|\lesssim 
\|\varphi_n\|_\Lone
+\|\varphi^+_n\|_\Lone+\|\varphi^-_n\|_\Lone\leq 1+2A$.
We then deduce, since $\varphi_n^+\leq 0$, that
$$\langle\mu,\varphi_n^-\rangle = -\langle\mu,\varphi_n\rangle + c_n +
\langle\mu,\varphi_n^+ \rangle \lesssim -n^3+1+2A.$$
It follows that $\langle \mu, \varphi \rangle =-\infty$ if
$\varphi:=\sum n^{-2} \varphi_n^-$. 
This is a contradiction because the series  $\sum n^{-2} \varphi_n^-$
converges to a q.p.s.h. function.
\end{proof}
\ 
\\
{\bf 2.2 Some properties of linear maps} 
\\
\ 
\\
Recall that a {\it Jordan block}
$J_{\lambda,m}$ is a square matrix
$(a_{i,j})_{1\leq i,j\leq m}$ such that
$a_{ij}=\lambda$ if $i=j$, $a_{ij}=1$ if $j=i+1$ and $a_{ij}=0$ 
otherwise. 
If $\lambda\not=0$, the entry of index
$(1,m)$ of $J_{\lambda,m}^n$ is 
equal to ${n\choose m-1}\lambda^{n-m+1}$. This is
the only entry of order
$n^{m-1}|\lambda|^n$,
the other ones have order at most 
 $n^{m-2}|\lambda|^n$.
We have 
$$\|J_{\lambda,m}^n\|\sim
{n \choose m-1}|\lambda|^{n-m+1} \sim n^{m-1}|\lambda|^n.$$ 
The eigenspace of $J_{\lambda,m}$
associated to the unique eigenvalue $\lambda$ is a complex line.
\par
If $E$, $E'$ are real (resp. complex) vector spaces, denote by
$\End(E,E')$ the space of $\R$-linear (resp. $\C$-linear)
maps from $E$ onto $E'$.
When $E$ and $E'$ are real vector spaces, we have
$$\End(E,E')\otimes_\R\C=\End(E\otimes_\R \C,E'\otimes_\R \C).$$
Hence, we can identify $\End(E,E')$ with a real vector subspace of $\End(E\otimes_\R
\C,E'\otimes_\R \C)$. 
\par
Consider a complex space 
$E$ and an invertible linear map $\Lambda\in \End(E,E)$. The space 
$E$ is the direct sum of the invariant complex subspaces
$E=E_1\oplus \cdots \oplus E_r$ with $\dim E_i=m_i$. The
restriction of
$\Lambda$ to $E_i$ is defined 
by a Jordan block
$J_{\lambda_i,m_i}$. We can assume that 
$(|\lambda_1|,m_1)\geq\cdots\geq (|\lambda_r|,m_r)$.
Let $F_i$ denote the eigenspace of
$\Lambda_{|E_i}$. It is a complex line. 
Let $E_i'$ be the hyperplane generated by the
first $(m_i-1)$ vectors of the basis of $E_i$ associated to the Jordan
form.
\par
We say that $J_{\lambda_i,m_i}$
is {\it dominant} if $(|\lambda_i|,m_i)= (|\lambda_1|,m_1)$.
In which case, we say that $\lambda_i$ 
is a {\it dominant eigenvalue} of $\Lambda$.
Assume that 
$J_{\lambda_1,m_1}$, $\ldots$, $J_{\lambda_\nu,m_\nu}$ are the dominant
Jordan blocks. 
It is clear that $\|\Lambda^n\|\sim
n^{m_1-1}|\lambda_1|^n$ and that for any vector $v\not\in
E_1'\oplus \cdots\oplus E_\nu'\oplus 
E_{\nu+1}\oplus\cdots \oplus E_r$, we have $\|\Lambda^n v\|\sim
n^{m_1-1}|\lambda_1|^n$.
The positive number 
$\lambda:=|\lambda_1|$ is the {\it spectral radius} of
$\Lambda$. The integer $m:=m_1$ is called the {\it multiplicity of
the spectral radius}. 
\par
We say that $F:=F_1\oplus \cdots
\oplus F_\nu$ is the {\it dominant eigenspace} and that
$F':=\oplus F_i$
with $1\leq i\leq \nu$, $\lambda_i=\lambda$, is the {\it 
strictly dominant eigenspace} of $\Lambda$.
These spaces are 
invariant under $\Lambda$. For any $1\leq j\leq \nu$, there is
a unique 
$\theta_j\in \mathbb{S}:=\R/2\pi\Z$ such that
$\lambda_j=\lambda\exp(i\theta_j)$. We say that 
$\theta:=(\theta_1,\ldots,\theta_\nu)\in \mathbb{S}^\nu$
is the {\it dominant direction} of $\Lambda$;
the dominant direction of $\Lambda^n$ is equal to $n\theta$. 
Denote by $\Theta$ the closed subgroup
of $\mathbb{S}^\nu$ generated by $\theta$.
It is a finite union
of real tori.
The orbit of each point
$\theta'\in \Theta$ under the translation
$\theta'\mapsto \theta'+\theta$ is dense in $\Theta$.  
If $\lambda_i=\lambda$ for every $1\leq i\leq\nu$,
we have $F=F'$, $\theta=0$ and $\Theta=\{0\}$.
\par
For every linear map
$L=(L_1,\ldots,L_r)$ in
$$\End(E,E)=\End(E,E_1)\oplus \cdots \oplus \End(E,E_r)$$
let
$$\exp(-in\theta)L:=\big(\exp(-in\theta_1)L_1,\ldots,
\exp(-in\theta_\nu)L_\nu,L_{\nu+1},\ldots,L_r\big).$$
We also define
$$\Lambda_n:=\frac{\Lambda^n}{n^{m-1}\lambda^n} \mbox{\ \ and \ \ }
\Lambda'_N:=\frac{1}{N}\sum_{n=1}^N \Lambda_i.$$
Let $\pi:F\longrightarrow F'$ denote the canonical projection.
The proof of the following proposition is left to the reader. 
One can deduce it from the proof of Proposition 2.4 if we take $K=E$ and 
$u=g=\id$. 
\begin{proposition} There exists a unique linear map
$\Lambda_\infty:E\longrightarrow F$ such that 
$$\|\exp(-in\theta)\Lambda_n-\Lambda_\infty\|=\O(1/n)
\ \mbox{ and  }\ \|\Lambda'_N-\pi\circ\Lambda_\infty\|=\O(\log N/N).$$
Moreover, $\pi\circ \Lambda_\infty$ is
of rank $\dim F'$. If $\Lambda$ preserves a convex cone $\K$ which
generates $E$ and satisfies
$-\overline \K\cap \overline \K =\{0\}$, then $\lambda$ is 
a dominant eigenvalue of $\Lambda$ with an eigenvector in $\overline \K$.
\end{proposition}
The last property of Proposition 2.3 is the Perron-Frobenius Theorem.
\begin{proposition} 
Let $K$ be a metric space, and let 
$u:K\longrightarrow E$ be a $\nu$-H\"older continuous vector-valued
function. Let $g:K\longrightarrow K$ be a Lipschitz map
such that $|g(x)-g(y)|\leq M |x-y|$ where $M>1$ is a constant. 
Assume that the spectral radius $\lambda$ of $\Lambda$ is strictly larger than $1$. 
Define
$$v_n:=\frac{1}{n^{m-1} \lambda^n}\sum_{j=1}^n \Lambda^j \circ u
\circ g^{n-j}
\mbox{ \ and \ }
w_N:=\frac{1}{N}\sum_{n=1}^N v_n.$$
Then there exist vector-valued functions $v$ and $w$ on $K$ 
such that 
$$\|\exp(-in\theta)v_n-v\|_\infty=\O(1/n)\ \mbox{ and }\ 
\|w_N-w\|_\infty=\O(\log N/N).$$
Moreover, $v$ and $w$ are $\nu'$-H\"older continuous for every
$\nu'>0$ with
$\nu' \leq \nu$ and $\nu'<\log \lambda/\log M$. 
\end{proposition}
\begin{proof} We can assume that the matrix $A$ of $\Lambda$ 
is the Jordan matrix $J_{\eta,m}$ with
$\eta=\exp(i\theta)\lambda$ and
$\theta\in\mathbb{R}/2\pi\Z$. Every entry of 
$A^n$ is of order at most 
$n^{m-2}\lambda^n$ except the one of index 
$(1,m)$ which is of order $n^{m-1}\lambda^n$. Since the functions
$u\circ g^n$ are uniformly bounded, every entry 
of $A^n$, whose order is smaller or equal to
$n^{m-2}\lambda^n$, does not contribute to $\lim v_n$ or $\lim w_N$.
This follows from the estimate $\sum_{j=0}^n j^{m-2}\lambda^j
\lesssim n^{m-2}\lambda^n$ for $\lambda>1$. Hence we can suppose that  
all the coordinate functions of $u$ are zero except the last one. 
\par
Let $u^+$ be the last coordinate function
of $u$. Define $u^+_j:=u^+\circ g^j$ for $j\geq 0$ and $s_j:= {j\choose m-1} 
\eta^{j-m+1}$ for $j\geq m-1$ with the convention that ${0\choose 0}=0$. 
The first coordinate functions of $v_n$ and $w_n$ are  
$$v^+_n:=\frac{1}{n^{m-1}\lambda^n}\sum_{j=m-1}^n
s_j u^+_{n-j}
\ \ \mbox{ and } \ \ 
w_N^+:=\frac{1}{N}\sum_{n=1}^N v_n^+.$$
It is sufficient to study the sequences of functions
$(v^+_n)$ and $(w_N^+)$.
\par
We show first that $\|v^+_{n+1}-\exp(i\theta)v^+_n\|_\infty\lesssim 
n^{-2}$. Observe that $\|u^+_j\|_\infty\leq\|u^+\|_\infty$ and
$$S_{j,n}:=\left|\frac{s_{j+1}}{(n+1)^{m-1}\lambda^{n+1}} - 
\frac{\eta s_j}{n^{m-1}\lambda^{n+1}}\right| \lesssim 
\frac{n-j+1}{n^2\lambda^{n-j}}.$$
Hence
\begin{eqnarray*}
\lefteqn{|v^+_{n+1}-\exp(i\theta)v^+_n|=}\\ 
& = & 
\left| \frac{1}{(n+1)^{m-1} \lambda^{n+1}}
\sum_{j=m-2}^n s_{j+1} u^+_{n-j}
- \frac{1}{n^{m-1}\lambda^{n+1}} \sum_{j=m-1}^n
\eta s_j u^+_{n-j} \right|\\
& \leq & \frac{|s_{m-1}|\|u^+_{n-m+2}\|_\infty}{(n+1)^{m-1}
\lambda^{n+1}}+
\sum_{j=m-1}^n S_{j,n} \|u_{n-j}^+\|_\infty 
 \lesssim  n^{-2}.
\end{eqnarray*}
Thus, the sequence $\exp(-in\theta) v_n^+$ converges 
uniformly to a function $v^+$ and $\|\exp(-in\theta)v^+_n-v^+\|_\infty
=\O(1/n)$. 
It follows that 
$\|w^+_N-w^+\|_\infty = \O(\log N /N)$ with $w^+=0$
if $\theta\not=0$ and $w^+=v^+$ 
if $\theta=0$.
\par
We only have to prove that
$v^+$ is H\"older continuous.
We will use the additive notation for the distance on $K$.
Let $x,y\in K$ and $\delta:=|x-y|$. First, consider the case where
$\nu=\log\lambda/\log M$, so $\lambda= M^\nu$.
Since $u^+$ is $\nu$-H\"older continuous, there exists
$c>0$, independent of $x$, $y$ such that
$$|u^+\circ g^j(x)-u^+\circ g^j(y)|\leq c|g^j(x)-g^j(y)|^\nu.$$
Hence
$$|u^+\circ g^j(x) -u^+\circ g^j(y)|\leq cM^{j\nu}\delta^\nu.$$
We also have $|u^+\circ g^j(x) -u^+\circ g^j(y)|\leq
2\|u^+\|_\infty$.
Let $q$ be the integer part of 
$-\log\delta/\log M$. We have the following estimates:
\begin{eqnarray*}
|v^+_n(x)-v^+_n(y)| & \leq &
\frac{1}{n^{m-1}\lambda^n}\sum_{j=m-1}^n
s_j |u^+\circ g^{n-j}(x)-u^+\circ g^{n-j}(y)|\\
& \lesssim& \frac{1}{n^{m-1}\lambda^n}\sum_{j=m-1}^n j^{m-1}\lambda^j 
|u^+\circ g^{n-j}(x)-u^+\circ
g^{n-j}(y)|\\
& \lesssim & \sum_{j=0}^\infty \lambda^{-j}|u^+\circ g^j(x) -u^+\circ
g^j(y)| \\
& \lesssim & \sum_{j=0}^q\lambda^{-j} M^{j\nu}\delta^\nu
+\sum_{j=q+1}^\infty \lambda^{-j}\\
& \lesssim & q\delta^\nu
+\lambda^{-q} \lesssim - (\log\delta)\delta^\nu.
\end{eqnarray*}
Consequently, $|v^+_n(x)-v^+_n(y)|\lesssim \delta^{\nu'}$. 
In the limit, we get $|v^+(x)-v^+(y)|\lesssim |x-y|^{\nu'}$.
This is the required inequality.
\par
The case where $\nu>\log\lambda/\log M$ is derived from the case
$\nu=\log\lambda/\log M$. If $\nu<\log\lambda/\log M$,
we have $M^\nu<\lambda$.
In the same way, we obtain
$$|v^+_n(x)-v^+_n(y)| \lesssim\sum_{j=0}^\infty \lambda^{-j}M^{j\nu}
\delta^\nu \lesssim \delta^\nu.$$
\end{proof}
\
\\
{\bf 2.3 Action on cohomology groups}
\\
\ 
\\
By Hodge theory \cite{GriffithsHarris, Demailly2}, 
we have the decomposition
$$\H^n(X,\C)=\sum_{p+q=n} \H^{p,q}(X,\C)$$
where $\H^{p,q}(X,\C)$ is the subspace of $\H^n(X,\C)$ spanned by classes
of closed $(p,q)$-forms. The group $\H^{p,q}(X,\C)$ is isomorphic to the 
Dolbeault cohomology group of bidegree $(p,q)$. We also have 
$\H^{p,q}(X,\C)=\overline{\H^{q,p}(X,\C)}$.
Define
$$\H^{p,p}(X,\R):=\H^{p,p}(X,\C)\cap\H^{2p}(X,\R).$$
Then
$$\H^{p,p}(X,\C)=\H^{p,p}(X,\R)\otimes_\R\C.$$
The {\it
K\"ahler cone} $\K_X$ of 
$\H^{1,1}(X,\R)$ is the cone of classes of K\"ahler forms
on $X$. In what follows, $\K_X^*$ will denote the cone of classes of positive closed
$(1,1)$-currents. These cones are convex;
$\K_X$ is open; $\K_X^*$ is
closed. We also have
$\overline \K_X\subset \K_X^*$ and $(-\K_X^*)\cap
\K_X^*=\{0\}$. 
\par
We will construct some invariant currents by
induction starting from an invariant current. It is necessary to introduce 
cohomology groups relative to a positive closed current: the manifold $X$ is 
replaced by a current $T$.  
Let $T\not=0$ be a positive closed current of bidegree $(s,s)$ with
$0\leq s\leq k-1$. 
For any $1\leq p\leq k-s$, let
$N^{p,p}(T,\R)$ denote the space of 
classes $[\alpha]\in \H^{p,p}(X,\R)$
which satisfy
$$[T]\wedge[\alpha]=0 \ \mbox{ in }\ \H^{p+s,p+s}(X,\R).$$
Let $N^{1,1}_\nu(T,\R)$ denote the space of classes $[S]$
where 
$S$ is a closed $(1,1)$-current (not necesarily positive)
with a $\nu$-H\"older continuous potential
such that $T\wedge S=0$. 
Define
$$\H^{p,p}(T,\R):=\frac{\H^{p,p}(X,\R)}{N^{p,p}(T,\R)}
\ \mbox{ and }\ 
\H_\nu^{1,1}(T,\R):=\frac{\H^{1,1}(X,\R)}{N^{1,1}_\nu(T,\R)}.$$
The map $\pi_T:\H^{p,p}(T,\R)\longrightarrow
\H^{p+s,p+s}(X,\R)$, given by $[\alpha]\mapsto [T]\wedge[\alpha]$,
is injective. Hence
$\|[T]\wedge.\|$ is a norm on
$\H^{p,p}(T,\R)$ if $\|.\|$ is a norm on
$\H^{p+s,p+s}(X,\R)$.
\par
Let $f$ be a holomorphic automorphism of $X$ such that 
$f^* T=\lambda_T T$ where $\lambda_T>0$ is a constant.
We define the map $f^*$ on $\H^{p,p}(X,\R)$ and $\H^{p,p}(T,\R)$ by
$f^*[\alpha]=[f^*(\alpha)]$ for every closed $(p,p)$-current
$\alpha$ on $X$. 
Let 
$$\lambda_{p,n}(T):=
\lambda_T^{-n}\|[(f^n)^*(T\wedge\omega^p)]\|=\|[T\wedge
(f^n)^*\omega^p]\|$$
and
$$\lambda_p(T):=\limsup_{n\rightarrow\infty}
\lambda_{p,n}(T)^{1/n}=\limsup_{n\rightarrow\infty}\left(\int_X 
T\wedge (f^n)^*\omega^p\wedge \omega^{k-s-p}\right)^{1/n}.$$
Observe that $\lambda_p(T)$ 
depends neither on the K\"ahler form $\omega$ nor on the norm
$\|. \|$ on $\H^{p+s,p+s}(X,\R)$. 
The following proposition follows from the discussion on Jordan 
forms for linear maps and from  Proposition 2.3 
(see also \cite{DinhSibony4,DinhSibony6}).
\begin{proposition} Let $X$, $f$ and $T$ be as above.
Then, $\lambda_p(T)$ is the spectral radius of $f^*$
on $\H^{p,p}(T,\R)$. If $l_p(T)$ is its multiplicity, then
$$\lambda_{p,n}(T)\sim n^{l_p(T)-1} \lambda_p(T)^n.$$
In particular, $\lambda_{p,n}(T)^{1/n}$
converge to $\lambda_p(T)$. 
\end{proposition}
\begin{proposition} Let $X$, $T$, $f$ be as above. Then for 
every $p_1\geq 1$ and $p_2\geq 1$ such that
$p_1+p_2\leq k-s$, we have $\lambda_{p_1+p_2}(T)\leq 
\lambda_{p_1}(T)\lambda_{p_2}(T)$.
In particular, 
$\lambda_1(T)^p\geq \lambda_p(T)$ for $1\leq p\leq k-s$
and $\lambda_1(T)^{k-s}\geq \lambda_T^{-1}$.
\end{proposition}
\begin{proof} 
Propositions 2.3 and 2.5 imply the existence of
$[\alpha_1]\in \H^{p_1,p_1}(T,\R)$, $[\alpha_2]\in \H^{p_2,p_2}(T,\R)$
such that
$$\frac{[(f^n)^*\omega^{p_1}]}{n^{l_{p_1}(T)-1}\lambda_{p_1}(T)^n}
\longrightarrow [\alpha_1]
\mbox{\ \ and \ \ }
\frac{[(f^n)^*\omega^{p_2}]}{n^{l_{p_2}(T)-1}\lambda_{p_1}(T)^n}
\longrightarrow [\alpha_2].$$
Hence
$$\frac{ [(f^n)^*\omega^{p_1+p_2}]}{n^{l_{p_1}(T)+
l_{p_2}(T)-2} \lambda_{p_1}(T)^n
\lambda_{p_2}(T)^n} \longrightarrow
[\alpha_1]\wedge [\alpha_2]$$
in $\H^{p_1+p_2,p_1+p_2}(T,\R)$.
On the other hand, there exists a non-zero class
$[\alpha]\in \H^{p_1+p_2,p_1+p_2}(T,\R)$
such that 
$$\frac{[(f^n)^*\omega^{p_1+p_2}]}{n^{l_{p_1+p_2}(T)-1} 
\lambda_{p_1+p_2}(T)^n}
\longrightarrow [\alpha].$$  
The property that $[\alpha]\not=0$ implies that
$\lambda_{p_1+p_2}(T)\leq\lambda_{p_1}(T)\lambda_{p_2}(T)$.
\par
The inequality $\lambda_1(T)^p\geq \lambda_p(T)$ is clear.
Since $f$ is an automorphism,
the mass of the measure $(f^n)^*(T\wedge \omega^{k-s})$ is equal to the mass 
of $T\wedge\omega^{k-s}$. Thus, $\lambda_T\lambda_{k-s}(T)=1$ and 
$\lambda_1(T)^{k-s}\geq \lambda_{k-s}(T)=\lambda_T^{-1}$.  
\end{proof}
\par
The space $\H_\nu^{1,1}(T,\R)$ is
invariant under $f^*$. Let $\rho_\nu(T)$
denote the spectral radius of $f^*$ on
$\H_\nu^{1,1}(T,\R)$ and
$m_\nu(T)$
its multiplicity. Since
$N^{1,1}_\nu(T,\R)\subset N^{1,1}(T,\R)$
we have
$$\big(\lambda_1(T),l_1(T)\big)\leq \big(\rho_\nu(T),m_\nu(T)\big).$$
We will prove later (see Lemma 3.3) that
if $\lambda_1(T)>1$, and if $\nu$ is small enough, then the last inequality
is in fact an equality.
\par
When $T$ is the integration current on $X$, we define
$$d_{p,n}:= \|[(f^n)^*\omega^p]\| \mbox{\ \ and\
\ }
d_p:=\lim_{n\rightarrow\infty}
\sqrt[n]{d_{p,n}}.$$
We also have
\begin{eqnarray}
d_p & = & \lim_{n\rightarrow \infty}
\left(\int_X (f^n)^*\omega^p\wedge \omega^{k-p}\right)^{1/n}.
\end{eqnarray}
The numbers $d_p$ are called the {\it dynamical degrees} of $f$. The
last one $d_t:=d_k$ is the {\it topological degree} of
$f$ which is equal to 1 because $f$ is an automorphism. 
It is noticed in \cite{Guedj} that an inequality of Khovanskii-Teissier-Gromov 
shows that $p\mapsto \log d_p$ is concave,
hence the sequence 
$(d_{p-1}/d_p)_{1\leq p\leq k}$ is increasing 
\cite{Khovanskii, Tessier, Gromov2}. 
In particular, there exist $m$, $m'$, $1\leq m\leq m'\leq k-1$, such that 
$$1\leq d_1<\cdots<d_m=\cdots =d_{m'}> \cdots >d_k=1.$$
\par
From Proposition 2.5, we know that $d_p$ is the spectral radius of $f^*$ acting on 
$\H^{p,p}(X,\R)$. Let $l_p$ denote its multiplicity.
Since $f^*$ preserves the cone of classes of positive closed $(p,p)$-currents, by
Proposition 2.3, $d_p$ is a dominant eigenvalue of $f^*$ on 
$\H^{p,p}(X,\R)$. The relation
$$\int_X (f^n)^*\omega^p\wedge\omega^{k-p}= 
\int_X \omega^p\wedge (f^n)_*\omega^{k-p}$$
implies that the spectral radius of $f_*$ on $\H^{k-p,k-p}(X,\R)$
is also equal to $d_p$ and its multiplicity is equal to $l_p$. 
\par
According to the Gromov-Yomdin theorem, the topological 
entropy $\h(f)$ of $f$
is equal to $\max_{1\leq p\leq k}\log d_p$ 
\cite{Gromov1, Yomdin, DinhSibony4, DinhSibony6}. In particular, if $\h(f)>0$,
we have $\max d_p>1$. It follows from 
Proposition 2.6 that $d_1>1$. It is shown in 
\cite{DinhSibony5} that the map which 
associates to a holomorphic endomorphism
of $X$ its topological entropy has discrete image in $[0,+\infty[$. 
It is well known that if $f$ belongs to the component of
the identity in the automorphism group of $X$, then $\h(f)=0$.
\section{Relative Green currents}
Let $f$ be a holomorphic automorphism of  a compact 
K\"ahler manifold  $(X,\omega)$ of dimension $k$.
Let $T$ be a positive closed $(s,s)$-current,
$0\leq s\leq k-1$, on $X$,
which satisfies a relation $f^*T=\lambda_T T$, $\lambda_T>0$. 
When $s=0$, $T$ is a multiple of the integration current
on the manifold $X$ and $\lambda_T=1$.
Define
$$M_n:=\|\mbox{D}f^n\|_\infty \mbox{ \ and \ }
M:=\lim_{n\rightarrow \infty} M_n^{1/n}$$
where $\mbox{D} f^n$ is the differential of $f^n$. The constant $M$ is 
independent 
of the metric on $X$. 
\par
Assume that $\lambda_1(T)>1$. Let $\Gamma(T)$ denote the cone
of $(s+1,s+1)$-currents $T\wedge S$ where $S$
is a closed $(1,1)$-current (not
necessarily positive) satisfying the following properties:
\begin{enumerate}
\item $S$ has a $\nu$-H\"older 
continuous potential for every $\nu$ such that
$0<\nu<\log\lambda_1(T)/\log M$;
\item $T\wedge S$ is a positive current;
\item $T\wedge f^*S=\lambda_1(T) T\wedge S$.
\end{enumerate}
We will denote by $\R\Gamma(T)$ the real space generated by $\Gamma(T)$.
We can now describe $\Gamma(T)$. 
\begin{theorem} Let $X$, $T$ and $f$ be as above.
Assume that $\lambda_1(T)>1$.
Then,  $\Gamma(T)$ is a closed finite dimensional 
cone with non zero elements.
Let $R$ be a closed real $(1,1)$-current 
with a continuous potential.
Then, the sequence of currents
$$\frac{1}{N}\sum_{n=1}^N \frac{T\wedge (f^n)^*R}{n^{l_1(T)-1}
\lambda_1(T)^n}$$
converges to a current in $\R\Gamma(T)$,
which depends only on the class $[R]$ in $\H^{1,1}(X,\R)$.
If $[R]$ belongs to $\K_X$,
the limit current belongs to $\Gamma(T)\setminus\{0\}$.
\end{theorem}
\par
Fix $\nu$ such that $0<\nu <\log\lambda_1(T)/\log M$. 
Replacing
$f$ by $f^n$ with $n>>0$, we can assume that
$0<\nu <\log\lambda_1(T)/\log M_1$.
Let  $m$ be the largest integer such that
$[\omega]$, $[f^*\omega]$, $\ldots$, $[(f^{m-1})^*\omega]$
are linearily independent in
$\H^{1,1}_\nu(T,\R)$.
Then there exist real numbers $a_0$, $\ldots$, $a_{m-1}$ such that
we have in $\H_\nu^{1,1}(T,\R)$ 
$$[(f^m)^*\omega] = a_{m-1}[(f^{m-1})^* \omega] +\cdots +
a_0[\omega].$$ 
Let $E$ be the subspace of $\H^{1,1}_\nu(T,\R)$ generated
by $[\omega]$, $\ldots$, $(f^{m-1})^*[\omega]$. 
These $m$ classes form 
a basis ${\cal B}$ of $E$ and $E$ is invariant under
$f^*$. Denote by $\Lambda$ the restriction of $f^*$ to $E$. 
The notations $\theta$, $\theta'$, $\Theta$, $\rho_\nu(T)$, $m_\nu(T)$ were
introduced in Section 2.
\par
The matrix of $\Lambda$ with respect to ${\cal B}$ is
$$A:=\left(
\begin{array}{ccccc}
0 & 0 & \cdots & 0 & a_0\\
1 & 0 & \cdots & 0 & a_1\\
0 & 1 & \cdots & 0 & a_2\\
\vdots & \vdots & \ddots & \vdots & \vdots\\
0 & 0 & \cdots & 1 & a_{m-1}
\end{array}
\right).
$$
Since $E$ contains a K\"ahler class, the spectral radius
of $\Lambda$ and of $A$ are equal to $\rho_\nu(T)
\geq \lambda_1(T)>1$. Moreover, their multiplicities are equal to
$m_\nu(T)$. We have
$$\|\Lambda^n\|=\|A^n\|\sim n^{m_\nu(T)-1}\rho_\nu(T)^n.$$ 
\begin{lemma} There exist a closed $(1,1)$-current 
$S$ and a continuous family of closed $(1,1)$-currents $S_{\theta'}$
with $\nu$-H\"older continuous potentials, for $\theta'\in\Theta$,
such that 
$T\wedge S\not =0$,
$T\wedge S_{\theta'} \not =0$, 
$T\wedge f^*S=\rho_\nu(T)T\wedge S$ and
$T\wedge f^*S_{\theta'}=\rho_\nu(T)T\wedge S_{\theta'+\theta}$.
Moreover, the sequences of positive closed currents 
$$\widetilde Z_N:=\frac{1}{N} \sum_{n=1}^N 
\frac{T\wedge (f^n)^*\omega}{n^{m_\nu(T)-1}\rho_\nu(T)^n}
\ \ \mbox{ \ and \ }\ \ 
Z_{n_i}:=\frac{T\wedge (f^{n_i})^*\omega}
{n_i^{m_\nu(T)-1}\rho_\nu(T)^{n_i}}$$
converge to $T\wedge S$ and to $T\wedge S_{\theta'}$ when $N\rightarrow
\infty$ and $n_i\rightarrow\infty$ with respectively $n_i\theta\rightarrow\theta'$. In
particular, if $\Theta$ is reduced to one point, $Z_n$
converge to $T\wedge S$.
\end{lemma}
\begin{proof}
From the definition of $\H^{1,1}_\nu(T,\R)$,
there exists a $(1,1)$-current 
$R$ with a $\nu$-H\"older continuous potential such that
$T\wedge R=0$ and 
$$\left[(f^m)^*\omega-\sum_{j=1}^m
a_{m-j} (f^{m-j})^*\omega\right]=[R] \ \mbox{ in }\ \H^{1,1}(X,\R).$$
Hence, there exists a $\nu$-H\"older continuous function $u$ such that
$$(f^m)^*\omega -\sum_{j=1}^m a_{m-j} (f^{m-j})^*\omega  =  R+\ddc u.$$
Then
$$T\wedge \left((f^m)^*\omega -\sum_{j=1}^m a_{m-j} (f^{m-j})^*\omega 
\right)  = T\wedge \ddc u.$$
Define
$$W_n:=\left(
\begin{array}{c}
(f^n)^*\omega\\
\vdots\\
(f^{n+m-2})^*\omega\\
(f^{n+m-1})^*\omega
\end{array}
\right)
\ \ \mbox{ and }\ \ 
U:=\left(
\begin{array}{c}
0\\
\vdots\\
0\\
u
\end{array}
\right).
$$
Then $W_{n+1}=f^*W_n$
and $T\wedge W_1=T\wedge BW_0+T\wedge \ddc U$ where $B$ is
the transpose of $A$. By induction,
we obtain
$$T\wedge W_n=T\wedge \left(B^n W_0 + \ddc
\sum_{j=1}^n B^{j-1} U\circ f^{n-j}\right).$$
Define
$$\overline W_n:=\frac{W_n}{n^{m_\nu(T)-1}\rho_\nu(T)^n}$$
and
$$V_n:=\frac{1}{n^{m_\nu(T)-1}\rho_\nu(T)^n} \left(B^n W_0 + \ddc
\sum_{j=1}^n B^{j-1} U\circ f^{n-j}\right).$$
Denote by $\overline W_n^+$ and $V_n^+$ the
first components of 
$\overline W_n$ and of $V_n$. 
We have $[\overline W_n^+]=[V_n^+]$
in $\H^{1,1}_\nu(T,\R)$.
Proposition 2.3 implies that the sequence
of classes $[\overline W_n^+]$ is bounded. 
Moreover, since $\omega$ is a K\"ahler form,
any cluster point of of this sequence is a non-zero class.
\par
Propositions 2.3 
and 2.4 imply that when $(n_i\theta)$ converges to $\theta'$, 
the sequence $V^+_{n_i}$ converges to a current $S_{\theta'}$
with a $\nu$-H\"older continuous potential. Moreover, $S_{\theta'}$ depends
on $\theta'$ but not on $(n_i)$.
The current $T\wedge S_{\theta'}$ is positive and closed.
Since $[S_{\theta'}]\not=0$ in $\H^{1,1}_\nu(T,\R)$,
from the definition of $\H^{1,1}_\nu(T,\R)$, 
we have $T\wedge S_{\theta'}\not=0$. 
It is clear
that $S_{\theta'}$
depends continuously on $\theta'\in \Theta$.
Since $\lim(n_i+1)\theta = \theta'+\theta$, we have
\begin{eqnarray*}
T\wedge f^*S_{\theta'} - \rho_\nu(T) T\wedge S_{\theta'+\theta}
& = & \lim \left(\frac{1}{\lambda_T}f^*Z_{n_i} - \rho_\nu(T) Z_{n_i+1}\right) \\
& = & \rho_\nu(T)
\lim \left(\left[\frac{n_i+1}{n_i}\right]^{m_\nu(T)-1}Z_{n_i+1}
- Z_{n_i+1}\right)  \\
& = &  0.
\end{eqnarray*}
Hence, $T\wedge f^*S_{\theta'} = \rho_\nu(T) T\wedge S_{\theta'+\theta}$.
\par
Propositions 2.3 and 2.4 imply that the sequence of currents 
$\frac{1}{N}\sum_{n=1}^N V^+_n$ converges to a current $S$
which has a $\nu$-H\"older continuous potential.
We also have $T\wedge S\not
=0$ because every limit value of 
$(T\wedge V^+_n)$ is a non-zero positive current. 
In the same way, we get 
$T\wedge f^*S=\rho_\nu(T) T\wedge S$.
\end{proof}
\begin{lemma} Under the hypothesis
of Theorem 3.1, we have
$$\big(\rho_\nu(T),m_\nu(T)\big)
=\big(\lambda_1(T),l_1(T)\big)$$ 
for every $\nu$ such that
$0<\nu<\log \lambda_1(T)/\log M$.
\end{lemma}
\begin{proof}
Recall that the map
$\pi:\H^{1,1}(T,\R)\longrightarrow \H^{s+1,s+1}(X,\R)$ defined
by  $\pi([\alpha]):=[T]\wedge[\alpha]$ is injective.
Consequently, the spectral radius
of $f^*$ on $\pi(\H^{1,1}(T,\R))$ is equal to
$\lambda_T\lambda_1(T)$ and its multiplicity is equal to
$l_1(T)$. We have seen that
the sequence of classes 
$$\frac{\big[(f^{n_i})^* (T\wedge \omega)\big]}
{n_i^{m_\nu(T)-1}\lambda_T^{n_i}\rho_\nu(T)^{n_i}}$$
converges to $[T\wedge S_{\theta'}]$ in $\H^{s+1,s+1}(X,\R)\setminus\{0\}$.
It follows from Proposition 2.5 that
$\big(\rho_\nu(T),m_\nu(T)\big)
=\big(\lambda_1(T),l_1(T)\big)$.
\end{proof}
\begin{lemma} With the assumptions of Theorem 3.1, 
let $R$ be a closed $(1,1)$-current
with a continuous potential such that
$[R]=\lambda [\omega]$ in
  $\H^{1,1}(X,\R)$ with $\lambda\in \R$. Then the sequence of
  currents
$$\frac{1}{N}\sum_{n=1}^N \frac{T\wedge (f^n)^*R}{n^{l_1(T)-1}
    \lambda_1(T)^n}$$
converges to $\lambda T\wedge S$. In particular, the
limit is $0$ if $[R]=0$ in $\H^{1,1}(X,\R)$.
\end{lemma}
\begin{proof}
Let $u$ be a continuous function
such that $R=\lambda\omega
  + \ddc u$. We have 
$$\frac{T\wedge (f^n)^*R}{n^{l_1(T)-1} \lambda_1(T)^n} = 
\lambda\frac{T\wedge (f^n)^*\omega}
{n^{l_1(T)-1} \lambda_1(T)^n} 
+T\wedge \ddc \left(\frac{u\circ f^n}{n^{l_1(T)-1}
  \lambda_1(T)^n}\right).$$
Since the $u\circ f^n$ are uniformly bounded, the last
relation implies that the sequence of currents
$$\frac{T\wedge (f^n)^*R}{n^{l_1(T)-1} \lambda_1(T)^n} -
  \lambda\frac{T\wedge (f^n)^*\omega}{n^{l_1(T)-1} \lambda_1(T)^n}$$
converges to $0$. It suffices to apply Lemmas 3.2 and 3.3.
The proof is valid for $R$ positive closed with a bounded potential.
\end{proof}
{\bf End of the proof of Theorem 3.1.} 
Let $\omega_1$, $\ldots$, $\omega_j$ be K\"ahler forms
such that the classes $[\omega_1]$, $\ldots$, $[\omega_j]$ generate
$\H^{1,1}(X,\R)$.
We can apply 
Lemma 3.4 to the forms $\omega_i$. 
This implies the convergence in Theorem 3.1.
\par
We now show that 
$\dim \Gamma(T)\leq\dim\H^{1,1}(X,\R)$. 
Otherwise, there exists a non-zero current
$T\wedge S$ in $\R\Gamma(T)$ with $[S]=0$. Let $u$ be a
continuous function such that $S=\ddc u$. We have 
$$T\wedge S= \lim \frac{T\wedge (f^n)^*S}{\lambda_1(T)^n} 
= \lim T\wedge \ddc\left(\frac{u\circ f^n}{\lambda_1(T)^n}\right)=0.$$ 
This is impossible. 
\par
Let $T^{(1)}$, $\ldots$, $T^{(r)}$ be a maximal linearily
independent set in $\Gamma(T)$. The cone $\Gamma(T)$ is equal to the
intersection of the space generated by the $T^{(i)}$ and the cone of 
positive closed $(s+1,s+1)$-currents. Hence,
it is closed.
\hfill $\square$
\\
\par
When $T$ is the integration current on $X$, 
Propositions 2.3, 2.5, 2.6 and Theorem 3.1 imply the following result
 (see \cite{Cantat, Guedj1} for the case of 
surfaces and \cite{BedfordLyubichSmillie4, 
FornaessSibony2, Sibony2} for polynomial automorphisms).
\begin{corollary} Let $(X,\omega)$ be a compact K\"ahler manifold of 
dimension $k$.
Let $f$ be a holomorphic automorphism of $X$, 
of positive topological entropy.
Then $d_1>1$ and there exists
a positive closed $(1,1)$-current $T_1$ satisfying 
$f^*T_1=d_1T_1$. Moreover, the potential of $T_1$ is
H\"older continuous and 
the class $[T_1]$ belongs to $\overline \K_X$.
\end{corollary}
\begin{remarks} \rm
The construction of invariant currents in Theorem 3.1 is still valid if we
restrict to an invariant subspace  
$E$ of $\H^{1,1}(X,\R)$.
We have to assume that the spectral radius of $f^*$ on the projection
$E'$ of $E$ in $\H^{1,1}(T,\R)$ is strictly larger than $1$.
The construction gives invariant currents with 
H\"older continuous potentials which are not necessarily positive;
there is a non-zero current if $f^*_{|E'}$ has a dominant real eigenvalue.
In particular, every invariant $(1,1)$-current $T$ with continuous potential such that
$f^*T=\lambda T$, $|\lambda|>1$, has a H\"older continuous potential. 
Let $\K_X^b$ be the cone of non-zero classes of 
positive closed $(1,1)$-currents
with bounded potentials. If $E\cap\K^b_X\not= \{0\}$, we obtain a positive
current with a H\"older continuous potential. 
\end{remarks}
\section{Green currents}
Let $f$ be as in Section 3.
Recall that $d_s$ is the spectral radius of $f^*$ on $\H^{s,s}(X,\R)$. 
Let $l_s$ denote its multiplicity. 
We consider $\Gamma_s$ the cone of PB positive closed $(s,s)$-currents
$T$ such that $f^* T = d_s T$. Let
$\R\Gamma_s$ denote the real space generated by $\Gamma_s$.
For every positive closed $(s,s)$-current $T$, define 
$$\Ccal(T):=\big\{S\mbox{ positive closed } (s,s)\mbox{-current, }
S\leq cT \mbox{ for some } c>0\big\}$$
and $\R\Ccal(T)$ the real space generated by $\Ccal(T)$.
Let $[\ ]$ denote the map which associates to a closed 
$(s,s)$-current its cohomology class in $\H^{s,s}(X,\C)$.
\begin{theorem} Assume that $d_s>d_{s-1}$.
Let $S$ be a PB closed real $(s,s)$-current.  
Then, the sequence of currents
$$S_N:=\frac{1}{N}\sum_{n=1}^N \frac{(f^n)^*S}{n^{l_s-1} d_s^n}$$
converges to a current in $\R\Gamma_s$
which depends only on the class $[S]$ in $\H^{s,s}(X,\R)$.
Moreover,
$\Gamma_s \not = \{0\}$ and 
the restriction of $[\ ]$ 
to $\R\Gamma_s$ is injective.
Every current in $\R\Gamma_s$ is PC.
If $T$ belongs to $\Gamma_s$, then 
the restriction of $[\ ]$ to $\R\Ccal(T)$ is injective. 
In particular, the cones $\Gamma_s$ and 
$\Ccal(T)$ are finite dimensional and closed. If $[T]$ is 
extremal in the cone of classes of positive closed $(s,s)$-currents then $T$ is
extremal in the cone of positive closed $(s,s)$-currents.
\end{theorem}
\begin{proof} Let $\Phi$ be a DSH $(k-s,k-s)$-current such that 
$\ddc\Phi=\Omega^+-\Omega^-$ where $\Omega^\pm$ are positive closed 
$(k-s+1,k-s+1)$-currents.
Then $[\Omega^+]=[\Omega^-]$ and $\|\Omega^+\|=\|\Omega^-\|$.
Assume that $\|\Omega^+\|=\|\Omega^-\|\leq 1$.
Proposition 2.1 implies the existence of a $(k-s,k-s)$-form 
$\Phi_0=\Phi_0^+-\Phi_0^-$ such that $\Phi_0^\pm\leq 0$,
$\ddc\Phi_0=\ddc\Phi$ and $\|\Phi_0^\pm\|_\DSH\leq A$.
The current $\Psi_0:=\Phi-\Phi_0$ is $\ddc$-closed.
Define $\Omega^\pm_n := (f^n)_* \Omega^\pm$. 
Recall that the spectral radius of $f_*$ on $\H^{k-s+1,k-s+1}(X,\C)$
is equal to $d_{s-1}$ and that its multiplicity is equal to $l_{s-1}$ 
(see Section 2). 
Fix $\epsilon$, $0<\epsilon<d_s-d_{s-1}$. We have
$\|\Omega^\pm_n\|\lesssim (d_s-\epsilon)^n$. 

Proposition 2.1 implies the existence
of $(k-s,k-s)$-forms $\Phi_n=\Phi_n^+-\Phi_n^-$ such that 
 $\ddc\Phi_n=\Omega^+_n-\Omega^-_n$, $\Phi_n^\pm\leq 0$ 
and $\|\Phi_n^\pm\|_\DSH\lesssim 
(d_s-\epsilon)^n$. 
If $S$ is PB and $\Phi$ is smooth or if $S$ is smooth and $\Phi$ is DSH as 
above, 
we have $|\langle S, \Phi_n^\pm\rangle|\lesssim 
(d_s-\epsilon)^n$ for every $n\geq 0$.
We define by induction the $\ddc$-closed form $\Psi_n$
as $\Psi_n:=f_*\Phi_{n-1}-\Phi_n$. They satisfy
$\|\Psi_n\|_\Lone
\lesssim (d_s-\epsilon)^n$ for $n\geq 1$. 
On the other hand, we have
$$(f^n)_*\Phi  =  (f^n)_*\Psi_0 + (f^n)_*\Phi_0
=   (f^n)_*\Psi_0 + (f^{n-1})_*\Psi_1
+ (f^{n-1})_*\Phi_1.$$
So by induction, we get
$$(f^n)_*\Phi =   (f^n)_*\Psi_0
+ \cdots +f_*\Psi_{n-1} +  \Psi_n  +  \Phi_n.$$

Since $X$ is K\"ahler, every closed form which is 
$\d$-exact is $\ddc$-exact \cite[p.41]{Demailly2}. 
Hence the $\ddc$-closed form $\Psi_n$
defines a linear form
on $\H^{s,s}(X,\R)$ by $[\alpha]\mapsto\int \Psi_n\wedge \alpha$ for every 
real closed $(s,s)$-form $\alpha$. The Poincar\'e duality allows to associate
to $\Psi_n$ a unique class 
$c_n$ in $\H^{k-s,k-s}(X,\R)$. 
For $n\geq 1$ we have
$$\|c_n\|\lesssim \|\Psi_n\|_\Lone\lesssim 
(d_s-\epsilon)^n.$$ 
Define 
$$b_n:= (f^n)_*c_0 + (f^{n-1})_*c_1 +\cdots + c_n \ \mbox{ and }\ 
B_N:= \frac{1}{N}\sum_{n=1}^N 
\frac{b_n}{n^{l_s-1} d_s^n}.$$
As in the proof of the Proposition 2.4, we can check 
that the sequence $(B_N)$ converges to a class 
$B\in \H^{k-s,k-s}(X,\C)$ such that 
$\|B\|\leq c\|\Phi\|_\DSH$ where $c>0$ is a constant.

Now, assume that $S$ is smooth and $\Phi$ is DSH. 
Since $S$ is closed, then
\begin{eqnarray}
\langle (f^n)^*S, \Phi\rangle  = \langle S, (f^n)_*\Phi\rangle 
=   \int [S]\wedge b_n
+ \langle S, \Phi_n\rangle
\end{eqnarray}
and
$$\langle S_N,\Phi\rangle = \int [S]\wedge B_N  + 
\frac{1}{N}\sum_{n=1}^N\frac{\langle S , \Phi_n\rangle}{n^{l_s-1} d_s^n}.$$
The second term in the right hand side of the last equality 
tends to zero because $\|\Phi_n\|_\DSH\lesssim (d_s-\epsilon)^n$.
Hence 
$$\lim \langle S_N,\Phi\rangle =\int [S ]\wedge B\leq c\|\Phi\|_\DSH.$$
It follows that $(S_N)$ converges to a PB current $S_\infty$
which depends only on the class $[S]$. 
It is clear that $f^*S_\infty=d_s S_\infty$. 
Hence $S_\infty$ belongs to $\R\Gamma_s$ (we can write $S$ and $S_\infty$ as 
differences of positive closed currents). 
Observe that if $S$ is strictly positive, by definition of $d_s$ and 
$l_s$, we have 
$[S_\infty]\not=0$. Hence $S_\infty$ 
is a non-zero positive current and $\Gamma_s\not=\{0\}$.
\par
Now assume that $S$ is PB (not necessarily smooth) and $\Phi$ is smooth. 
Then, $\Phi_n$ is continuous. If $S'$ is a smooth real $(s,s)$-form cohomologous to
$S$, we have
\begin{eqnarray*}
\langle S_N-S'_N,\Phi\rangle & = &  \int [S-S']\wedge B_N  + 
\frac{1}{N}\sum_{n=1}^N\frac{\langle S-S' , \Phi_n\rangle}{n^{l_s-1} d_s^n}\\
& = & \frac{1}{N}\sum_{n=1}^N\frac{\langle S-S' , 
\Phi_n\rangle}{n^{l_s-1} d_s^n}
\end{eqnarray*}
The last term tends to zero because $S-S'$ is PB and 
$\|\Phi_n\|_\DSH\lesssim (d_s-\epsilon)^n$.
It follows that $(S_N)$  converges to a PB current in $\R\Gamma_s$.
\par
Let $R\in \R\Gamma_s$ be a current such that $[R]=0$. 
Then, using identity (3) we get
$$|\langle R,\Phi \rangle |= d_s^{-n} |\langle (f^n)^* R, \Phi \rangle|
= d_s^{-n} |\langle R, \Phi_n \rangle| 
\lesssim d_s^{-n}(d_s-\epsilon)^n.$$
Therefore, $\langle R, \Phi\rangle = 0$ and hence $R=0$. It follows 
that the restriction of $[\ ]$ to $\R\Gamma_s$ is injective.
\par
Let $R\in \R\Gamma_s$ and $\Phi$ smooth. Using the identity 
$(f^n)^*[R]=d_s^n[R]$, we get 
\begin{eqnarray*}
\langle R,\Phi \rangle & = & d_s^{-n} \langle R, (f^n)_* \Phi \rangle \\
& = & \int [R]\wedge (c_0+d_s^{-1}c_1+\cdots+d_s^{-n}c_n) +  
d_s^{-n} \langle R,\Phi_n\rangle.
\end{eqnarray*}
Since $R$ is PB, when $n\rightarrow \infty$, we get
$\langle R,\Phi \rangle=\int [R]\wedge c_\Phi$ with 
$c_\Phi:=\sum_{n\geq 0} d_s^{-n} c_n$.
Following Proposition 2.1, $c_\Phi$ depends continuously
on $\Phi$. 
Hence, we can extend $R$ to a continuous linear form on $\Phi\in\DSH^{k-s}(X)$ by
$$\langle R,\Phi\rangle:=c_\Phi.$$
Hence $R$ is PC.
\par
We show that the restriction of $[\ ]$ to $\R\Ccal(T)$ is injective.
Let $R\in \R\Ccal(T)$ be a current such that $[R]=0$.
We have to prove that $R=0$. 
We can write $R=R^+-R^-$ with $R^\pm$ positive closed currents such that
$R^\pm\leq cT$ for a constant $c>0$.
Define $R_n^\pm:=d_s^n (f^n)_*R^\pm$ and $R_n:=R_n^+-R_n^-$. We have
$R_n^\pm\leq  cd_s^n(f^n)_*T=cT$.
For a smooth test form $\Phi$, we have
$|\langle T,\Phi_n^\pm\rangle|\lesssim (d_s-\epsilon)^n$. The domination
of $R_n^\pm$ and the negativity of $\Phi^\pm_n$ imply that 
$|\langle R_n,\Phi^\pm_n\rangle| \lesssim (d_s-\epsilon)^n$.
Since $[R_n]=0$, we obtain from (3) that
$$|\langle R,\Phi\rangle| = d_s^{-n}|\langle (f^n)^* R_n,\Phi\rangle|
 = d_s^{-n}|\langle R_n,\Phi_n\rangle|
\lesssim d_s^{-n}(d_s-\epsilon)^n.$$
Hence $\langle R,\Phi\rangle =0$ and $R=0$. 
This completes the proof of Theorem 4.1.
\end{proof}

\begin{corollary} Let $f$ be a holomorphic automorphism of a 
compact K\"ahler manifold $X$ of dimension
$k$. Assume that the
dynamical degrees of $f$ are all distinct. Then for every $s$, $1\leq s\leq k$,
there exists a non-zero PC positive closed $(s,s)$-current $T_s$ such that
$f^* T_s =c_s T_s$ with $c_s>0$. 
\end{corollary}
\begin{proof} The hypothesis implies, thanks to the 
Khovanskii-Tessier-Gromov convexity
theorem \cite{Khovanskii, Tessier, Gromov2} (see Section 2.3), the existence of
$m$ such that 
$$1<d_1<d_2<\cdots < d_m>d_{m+1}>\cdots
>d_k=1.$$ 
Using Theorem 4.1, 
we construct the current $T_s$ such that $f^* T_s =d_s T_s$  
for $1\leq s\leq m$.
The current $T_1$ can be constructed as in Corollary 3.5.
We now construct the other currents by induction 
using Theorem 3.1 for $f^{-1}$.
We construct
$(1,1)$-currents $S_i$, $1\leq i\leq k-m$, 
with H\"older continuous potentials and 
invariant currents $T_s$, $m+1\leq s\leq k$, 
of the form
$T_s=T_m\wedge S_1\wedge\ldots\wedge S_{s-m}$.
These currents
satisfy  $f^*T_s=c_s T_s$, $c_s>0$. 
Since $f$ is an automorphism, we necesarily have $c_k=1$. 
\par
In order to apply inductively Theorem 
3.1 for $f^{-1}$, we need only to verify 
that the first dynamical degree 
$\lambda_1(T_s)$ of 
$f^{-1}$, relative to $T_s$, is strictly larger 
than 1 for $m\leq s\leq k-1$. Following the last 
inequality of Proposition 2.6, it is sufficient to prove that $c_s>1$ 
for $m\leq s\leq k-1$. We have 
for every $\epsilon>0$
\begin{eqnarray*}
c_s^{-n} & \lesssim & \int (f^n)_*T_s\wedge \omega^{k-s}
=d_m^{-n}\int T_m\wedge (f^n)_*(S_1\wedge\ldots\wedge S_{s-m})\wedge 
\omega^{k-s}\\
& = & d_m^{-n}\int [S_1]\wedge\ldots\wedge [S_{s-m}] \wedge (f^n)^*[T_m\wedge 
\omega^{k-s}] \lesssim d_m^{-n} (d_{k-s+m}+\epsilon)^n.
\end{eqnarray*}
It follows that $c_s>1$ for $m\leq s\leq k-1$. This completes 
the induction step.
\par
One can check that the wedge product of a PC positive closed current with a 
current of bidegree $(1,1)$ with continuous potential is always PC.
\end{proof}  
\section{Mixing of the equilibrium measure}
In this section, using the methods developed above, we can construct, 
for automorphisms with distinct dynamical degrees, an equilibrium measure which is 
PC and mixing. We get the following result. 
\begin{theorem} Let $f$ be a holomorphic automorphism of a
compact K\"ahler manifold $X$ of dimension $k$. Assume that 
the dynamical degrees of $f$ are all distinct. 
Then $f$ admits a mixing PC invariant measure $\mu$.
Moreover, $\mu$ gives no mass to sets with small Hausdorff dimension. 
\end{theorem}
\par
We need the following variation of Ahlfors's estimate 
(see \cite{BedfordSmillie,Sibony2}).
\begin{lemma} Let $f$ be a holomorphic automorphism of $X$. Let
$T$ be a positive closed $(s,s)$-current such that
$f^*T=\lambda_T T$ with 
  $\lambda_T>0$. Assume that $\lambda_1(T)<1$. Then for every
smooth function $\psi\geq 0$, the limit values
of the sequence
$S_n:=\lambda_T^{-n}(f^n)^*(\psi T)$ are positive closed currents.
Moreover,
$\|\d S_n\|\rightarrow 0$ and $\|\ddc S_n\|\rightarrow 0$. If
  $S_{n_i}\rightharpoonup S$ and if $\sigma$
is a closed $(1,1)$-current with a continuous potential, then
$S_{n_i}\wedge \sigma \rightharpoonup S\wedge \sigma$.
\end{lemma}
\begin{proof} Let $\theta$ be a continuous
$(0,1)$-form. The
Cauchy-Schwarz inequality implies that
\begin{eqnarray*}
A_n & := & \left|\int (f^n)^* (\partial\psi) \wedge
  T \wedge \theta \wedge \omega^{k-s-1} \right|  \\
& \leq & \left|\int (f^n)^* (\partial\psi \wedge
  \overline{\partial\psi}) \wedge
  T \wedge \omega^{k-s-1} \right|^{1/2}
\left|\int \theta \wedge
  \overline \theta \wedge
  T \wedge \omega^{k-s-1} \right|^{1/2}\\
&\leq &
c \left|\int (f^n)^*\omega \wedge
  T \wedge \omega^{k-s-1} \right|^{1/2}
\left|\int \omega \wedge T \wedge \omega^{k-s-1} \right|^{1/2}
\end{eqnarray*}
if $i\partial \psi\wedge \overline{\partial\psi}$ and
$i\theta\wedge\overline\theta$ are bounded by $c\omega$, $c>0$.
It follows that $A_n\lesssim (\lambda_1(T)+\epsilon)^{n/2}$.
Since $\lambda_1(T)<1$, we have $\lim A_n=0$.
As a consequence,
$\lim\|\partial S_n\|=0$,
hence $\lim\|\d S_n\|=0$.
\par
To estimate $\|\ddc S_n\|$, one has just to observe that for $c>0$ large enough
$$-c(f^n)^*\omega \wedge T\leq \ddc (f^n)^*\psi \wedge T 
\leq c (f^n)^*\omega\wedge T.$$
\par
Let $u$ be a local continuous potential of $\sigma$
and $\theta$ be a test form.
Define $\psi_n:=(f^n)^*\psi$.
For the last assertion of this lemma, we have 
\begin{eqnarray*}
\langle \psi_n T\wedge \ddc u,\theta\rangle & = & 
\langle \ddc (uT),\psi_n\theta \rangle \\
& = & \langle \ddc (\psi_n T), u\theta\rangle + 
\langle \d(\psi_n T), u\dc\theta\rangle -\\
& & -\langle \dc (\psi_n T), u\d\theta\rangle
+ \langle \ddc (u\psi_n T), \theta\rangle 
\end{eqnarray*}
The first three terms tend to zero. Hence, 
$S_{n_i}\wedge \sigma\rightharpoonup S\wedge \sigma$.
\end{proof}
{\bf Proof of Theorem 5.1.}
We construct invariant currents $T_s$ as in Corollary 4.2. We choose
an extremal current $T_m$  in $\Gamma_m$.
We can write $T_s=T_m\wedge S_1\wedge \ldots \wedge S_{s-m}$ for $s\geq m+1$, 
where $S_i$ are closed
$(1,1)$-currents with H\"older continuous potentials. 
Define $\mu:=T_k$.
Hence $\mu$ is PC.
The estimate of Hausdorff 
dimension of $\mu$ uses classical arguments \cite{DinhSibony1,Sibony2}. 
If the potentials of $S_j$, $1\leq j \leq k-m$, are
$\alpha_j$-H\"older continuous, then $\mu$ gives no mass to sets
whose Hausdorff dimension is smaller than $\alpha_1+\cdots
+\alpha_{k-m}$.
\par
We show first that $\mu$ is 
ergodic. Let $\psi\geq 0$ be a smooth
test function.
Let $\tau$ be the limit of a sequence of currents
$n_i^{-1}\sum_{j=1}^{n_i}(\psi\circ f^j)T_m$. It is clear that
$f^*\tau=d_m\tau$ and
$\tau\leq \|\psi\|_\infty T_m$. 
Lemma 5.2 implies that $\tau$ is closed. Hence $\tau\in\Gamma_m$.
Since
$T_m$ is extremal in $\Gamma_m$, we have $\tau=cT_m$ for a constant $c$. 
We can now apply inductively Lemma 5.2.
Since the currents
$S_r$ have  continuous potentials and 
 $\lambda_1(T_s)<1$ for $m\leq s\leq k-1$,
$n_i^{-1}\sum_{j=1}^{n_i}(\psi\circ f^j)\mu$ converge
to $c\mu$. 
The invariance property of $\mu$ implies that $c=\|\mu\|^{-1}
\int \psi \d\mu$. This constant does not depend on the sequence 
$(n_i)$. Consequently, $n^{-1}\sum_{j=1}^n(\psi\circ f^j)\mu$ 
converge to $c\mu$. 
Hence $\mu$ is ergodic.
\par
We now prove that $\mu$ is mixing, which means 
$(f^n)^*\psi\mu\rightarrow c\mu$, $c=\|\mu\|^{-1}\int\psi\d\mu$, 
for every smooth function $\psi$.
Let $M$ denote the set of measures which are limite values
of the sequence $(f^{n*}\psi)\mu$ for some smooth function $\psi$.
Since ${\cal C}(T_m)$ is finite dimensional, Lemma 5.2 implies that 
$M$ is a finite dimensional space which contains $\mu$
and which is invariant under $f^*$ and $f_*$. 
Let $E$ denote the space of functions $\varphi\in\Ltwo(\mu)$ such that 
$\int\varphi \d\mu'=0$ for every $\mu'\in M$ and $E^\perp$ its orthogonal.
Observe that these spaces are invariant under $f^*$, $f_*$ and
that we have $\dim E^\perp =\dim M$. 
Moreover, in $E$, every function can be approximated
by smooth ones.
\par
We show that $\dim E^\perp=1$.
Since $f^*$ and $f_*$ preserve the scalar product in $\Ltwo(\mu)$, every eigenvalue 
of $f^*$ or $f_*$ has modulus equal to 1. Let $\varphi$ be an eigenvector
of $f_*$ associated to an eigenvalue $\lambda$. We have
$f_*|\varphi|=|\varphi|$. 
The ergodicity of $\mu$ implies that $|\varphi|$ is constant.
In particular, $\varphi^n\in\Ltwo(\mu)$ for every $n\geq 1$ and we have
$f_*\varphi^n=\lambda^n\varphi^n$. 
We claim that $E$ does not contain any eigenvector. Otherwise, there is 
a function $\varphi\in E\setminus\{0\}$ such that $f_*\varphi=\lambda\varphi$ with a 
$\lambda$ such that $|\lambda|=1$. We have for every smooth function $\psi$:
$$|\langle \varphi\mu,\psi\rangle|=|\langle (f^n)_*\varphi\mu,\psi\rangle|
=|\langle (f^{n*}\psi)\mu,\varphi\rangle|\rightarrow 0.$$
The last relation follows from the definition of $E$. We need of course to approach 
$\varphi$ by smooth functions in $E$. We get that 
$\varphi\mu=0$, hence $\varphi=0$. A contradiction.
\par
Let $\varphi$ be an eigenvector of $f_*$ in $E^\perp$ associated to an eigenvalue
$\lambda$. Then, $\varphi^n$ belongs to $E^\perp$ and is an eigenvector associated 
to $\lambda^n$ for every $n\geq 1$. Since $\dim E^\perp$ is finite, $\lambda$ 
is a root of unity. We have $f^{n*}\varphi=\varphi$ for some $n\geq 1$.
Since $\mu$ is ergodic, $\varphi$ is constant. Hence, $\lambda=1$.
Since  $f^*_{|E^\perp}$ preserves 
the scalar product, 
$\dim E^\perp=1$ and
$M$ is generated only by $\mu$. 
If $(f^{n_i*}\psi)\mu \rightarrow
c\mu$, we have $c=\|\mu\|^{-1}\int\psi\d\mu$. This constant 
does not depend on $(n_i)$. Hence $(f^{n*}\psi)\mu\rightarrow c\mu$ 
and $\mu$ is mixing.
\hfill $\square$

\small
Tien-Cuong Dinh and Nessim Sibony,\\
Math\'ematique - B\^at. 425, UMR 8628, 
Universit\'e Paris-Sud, 91405 Orsay, France. \\
E-mails: TienCuong.Dinh@math.u-psud.fr and
Nessim.Sibony@math.u-psud.fr.

\begin{thebibliography}{11}
%
%
\bibitem{BedfordLyubichSmillie4}
\textit{E. Bedford, M. Lyubich and J. Smillie}, Polynomial diffeomorphisms of
$\C^2$ IV: 
The measure of maximal entropy and laminar currents,
\textit{Invent.Math.}, \textbf{112} (1993), 77-125.
%
\bibitem{BedfordSmillie}
\textit{E. Bedford and J. Smillie}, Polynomial diffeomorphisms of
$\C^2$ III: 
Ergodicity, exponents and entropy of the equilibrium
measure, \textit{Math. Ann.}, \textbf{294} (1992), 395-420.
%
\bibitem{Blanchard}
\textit{A. Blanchard},
Sur les vari\'et\'es analytiques complexes, 
\textit{Ann. Sci. Ecole Norm. Sup.} (3),
\textbf{73} (1956), 157--202.
%
\bibitem{BostGilletSoule}
\textit{J.-B. Bost, H. Gillet, C. Soul\'e},
Heights of projective varieties and positive Green forms,
\textit{J. Amer. Math. Soc.}, \textbf{7} (1994), no. 4, 903--1027.
%
\bibitem{BriendDuval2}
\textit{J.Y. Briend et J. Duval}, Deux caract\'erisations de la mesure
d'\'equilibre d'un endomorphisme de
$\P^k(\C)$, \textit{IHES Publ. Math.}, \textbf{93} (2001), 145--159. 
%
\bibitem{Cantat}
\textit{S. Cantat}, Dynamique des automorphismes des surfaces K3,  
\textit{Acta Math.}, \textbf{187} (2001), no. 1, 1-57.
%
%
\bibitem{ClozelUllmo}
\textit{L. Clozel et E. Ullmo}, Correspondances modulaires et
mesures invariantes, \textit{J. reine angew. Math.}, \textbf{558} (2003), 47-83.
%
%
\bibitem{Demailly1}
\textit{J.P. Demailly}, Monge-Amp\`ere Operators, Lelong numbers and
Intersection theory in Complex Analysis and Geometry, \textit{Plemum
  Press} (1993), 115-193, \textit{(V. Ancona and A. Silva editors)}.
%
\bibitem{Demailly2}
\textit{J.P. Demailly}, Introduction \`a la th\'eorie de Hodge,
\textit{Panoramas et Synth\`eses}, \textbf{3} (1996), 1-111.
%
\bibitem{Demailly3}
\textit{J.P. Demailly}, Pseudoconvex-concave duality and regularization of
currents. Several complex variables (Berkeley, CA, 1995-1996), 233-271, 
\textit{Math. Sci. Res. Inst. Publ.}, \textbf{37}, Cambridge Univ. Press, 
Cambridge, 1999. 
%
\bibitem{Dinh}
\textit{T.C. Dinh}, Distribution des pr\'eimages et des points
p\'eriodiques d'une correspondance polynomiale, 
{\it Bull. Soc. Math. France}, to appear.
%
\bibitem{Dinh2}
\textit{T.C. Dinh}, Suites d'applications m\'eromorphes multivalu\'ees et 
courants laminaires, \textit{preprint}, 2003. arXiv:math.DS/0309421.
%
\bibitem{DinhSibony2}
\textit{T.C. Dinh et N. Sibony}, Dynamique des applications
d'allure polynomiale, {\it J. Math. Pures Appl.}, 
\textbf{82} (2003), 367-423. 
%
\bibitem{DinhSibony1}
\textit{T.C. Dinh et N. Sibony}, Dynamique des applications
polynomiales semi-r\'eguli\`eres, {\it Arkiv f\"or Mahematik}, {\bf 42} (2004),
61-85.
%
\bibitem{DinhSibony5}
\textit{T.C. Dinh et N. Sibony}, Groupes commutatifs
d'automorphismes d'une vari\'et\'e k\"ahl\'erienne compacte, 
\textit{Duke Math. J.}, {\bf 123} (2004), no. 2, 311-328. 
%
\bibitem{DinhSibony3}
\textit{T.C. Dinh et N. Sibony}, Distribution de valeurs d'une
suite de transformations m\'eromorphes et applications, {\it preprint},
2003. arXiv:math.DS/0306095.
%
\bibitem{DinhSibony4}
\textit{T.C. Dinh et N. Sibony}, Une borne sup\'erieure de l'entropie
topologique d'une application rationnelle, {\it Ann. of Math.}, to appear.
%
\bibitem{DinhSibony6}
\textit{T.C. Dinh et N. Sibony}, Regularization of currents and entropy, 
{\it Ann. Sci. Ecole Norm. Sup.}, to appear.
%
\bibitem{DinhSibony7}
\textit{T.C. Dinh et N. Sibony}, Dynamics of regular birational 
maps in $\P^k$, 
{\it J. Funct. Anal.}, to appear.
%
\bibitem{DinhSibony9}
\textit{T.C. Dinh et N. Sibony}, Decay of correlations and central limit 
theorem for meromorphic maps, {\it preprint}, 2004.
arXiv:math.DS/0410008.
%
\bibitem{FavreGuedj}
\textit{C. Favre et V. Guedj}, Dynamique des applications rationnelles
des espaces multiprojectifs, \textit{Indiana Univ.
  Math. J.}, \textbf{50} (2001), no. 2, 881-934.
%
\bibitem{FornaessSibony1}
\textit{J.E. Forn\ae ss and N. Sibony}, Complex H\'enon mappings in
$\C^k$ and Fatou-Biebebach domains, \textit{Duke Math. J.},
\textbf{65} (1992), 245-250.
%
\bibitem{FornaessSibony2}
\textit{J.E. Forn\ae ss and N. Sibony}, Complex dynamics in higher
dimension, \textit{in Complex potential theory}, (Montr\'eal, PQ,
1993), Nato ASI series Math. and Phys. Sci., vol. \textbf{C439},
Kluwer (1994), 131-186.
%
\bibitem{GilletSoule}
\textit{H. Gillet and C. Soul\'e}, Arithmetic intersection theory, 
\textit{I.H.E.S. Publ. Math.}, \textbf{72} (1990), 93--174 (1991).
%
\bibitem{GriffithsHarris} 
\textit{P. Griffiths and J. Harris}, Principles of algebraic geometry,
Wiley Classics Library, John Wiley \& Sons, Inc., New York, 1994. 
%
\bibitem{Gromov1}
\textit{M. Gromov}, On the entropy of holomorphic maps, 
\textit{Enseignement Math.}, \textbf{49} (2003), 217-235 
(\textit{manuscript}, 1977).
%
\bibitem{Gromov2}
\textit{M. Gromov}, Convex sets and K\"ahler manifolds,
\textit{Advances in differential geometry and topology}, 
Word Sci. Publishing, Teaneck, NJ, 1998, 1-38. 
%
\bibitem{Guedj1}
\textit{V. Guedj}, Dynamics of polynomial mappings of $\C^2$, 
\textit{Amer. J. Math.}, \textbf{124} (2002), no. 1, 75--106.
%
%
\bibitem{Guedj}
\textit{V. Guedj}, Ergodic properties of rational mappings with
large topological degree, \textit{Ann. of Math.}, to appear.
%
%
\bibitem{GuedjSibony}
\textit{V. Guedj et N. Sibony}, Dynamics of polynomial automorphisms
of $\C^k$, \textit{Arkiv f\"or Matematik}, \textbf{40} (2002), 207-243.
%
\bibitem{Khovanskii}
\textit{A. Khovanskii}, Fewnormials and Pfaff manifolds, {\it I.C.M. 1983}, 
{\it Warsawa} (1984), 549-565.
%
\bibitem{Margulis}
\textit{G. A. Margulis}, Discrete subgroups of semisimple Lie groups, 
Springer-Verlag, 1991.
%
\bibitem{Mazur}
\textit{B. Mazur}, The topology of rational points,
\textit{Expriment Math.}, \textbf{1} (1992), 35-45.
%
\bibitem{McMullen}
\textit{C.T. McMullen}, Dynamics on K3 surfaces: Salem numbers
and Siegel disks, \textit{J. reine angew. Math.}, \textbf{545}
(2002), 201-233.
%
\bibitem{Sibony2}
\textit{N. Sibony}, Dynamique des applications rationnelles de
$\mathbb{P}^k$, \textit{Panoramas et Synth\`eses}, {\bf 8} (1999), 97-185.
%
\bibitem{Tessier}
\textit{B. Tessier}, Bonnesen-type inequalities in algebraic geometry, 
\textit{Sem. on Diff. Geom.} (1982), 85-105. Princeton Univ. Press.
%
\bibitem{Voisin}
\textit{C. Voisin}, Intrinsic pseudovolume forms and
K-correspondences, \textit{preprint}, 2003. 
%
\bibitem{Yomdin}
\textit{Y. Yomdin}, Volume growth and entropy, \textit{Israel
  J. Math.}, \textbf{57} (1987), 285-300.
%
\end{thebibliography}
\end{document}